\documentclass[OJMO,NoBabel]{cedramhack}
\makeatletter
\renewcommand*{\thetable}{\arabic{table}}
\renewcommand*{\thefigure}{\arabic{figure}}
\let\c@table\c@figure
\makeatother

\vbadness=\maxdimen
\title[Trading off 1-norm against rank for ah-symmetric generalized inverses]{Trading off 1-norm and sparsity against rank for\\ linear models using mathematical optimization:\\
1-norm minimizing partially reflexive ah-symmetric\\ generalized inverses
}

\author[M. Fampa]{\firstname{Marcia} \lastname{Fampa}}
\address{Federal University of Rio de Janeiro}
\email{fampa@cos.ufrj.br}

\author[J. Lee]{\firstname{Jon} \lastname{Lee}}
\address{University of Michigan}
\email{jonxlee@umich.edu}

\author[G. Ponte]{\firstname{Gabriel} \lastname{Ponte}}
\address{Federal University of Rio de Janeiro}
\email{gabrielponte@poli.ufrj.br}

\thanks{M. Fampa was supported in part by CNPq grants 303898/2016-0 and 434683/2018-3.
J. Lee was supported in part by AFOSR grant FA9550-19-1-0175.
J. Lee and M. Fampa were supported in part by funding from the Simons
Foundation and the Centre de Recherches Math\'ematiques, through the Simons-CRM
scholar-in-residence program. G. Ponte was supported in part by CNPq PIBIC scholarship 149149/2020-4.}


\begin{abstract}
	The M-P (Moore-Penrose) pseudoinverse has as a key application the computation of  least-squares solutions of inconsistent  systems of linear equations. Irrespective of whether a given input matrix is sparse, its M-P pseudoinverse can be  dense, potentially leading to high computational burden, 
	especially when we are dealing with high-dimensional matrices. The  M-P pseudoinverse is uniquely characterized by four properties, but only two of them need to be satisfied for the computation of least-squares solutions. 
	Fampa and Lee (2018) and Xu, Fampa, Lee, and Ponte (2019)  propose local-search procedures to construct sparse block-structured generalized inverses that satisfy the two key M-P properties, plus one more (the so-called reflexive
	property).
	That additional M-P property is equivalent to imposing a minimum-rank condition on the generalized inverse.
	 (Vector) 1-norm minimization is used to induce sparsity and, importantly, to keep the magnitudes of entries under control for the generalized-inverses constructed.
		Here, we 
		investigate the trade-off between low 1-norm and low rank for generalized inverses that can be used in the computation of least-squares solutions. 
		We propose several algorithmic approaches that start from a $1$-norm minimizing generalized inverse that satisfies the two key M-P properties, and gradually decrease its rank, by iteratively imposing the reflexive property. The algorithms iterate until the generalized inverse has the least possible rank. During the iterations,  we produce intermediate solutions,
		trading off low 1-norm (and typically high sparsity) against low rank. 
\end{abstract}
\begin{document}

\maketitle

\section{Introduction}

The well-known M-P (Moore-Penrose) pseudoinverse
%, independently discovered  by  E.H. Moore and R. Penrose, 
is used in several linear-algebra applications, as for example, to compute least-squares solutions of inconsistent  systems of linear equations. If $A=U\Sigma V^\top$ is the real singular value decomposition of $A\in\mathbb{R}^{m\times n}$ (see \cite{GVL1996}, for example),
then the
M-P pseudoinverse of $A$ can be defined as $A^{\dagger}:=V\Sigma^{\dagger} U^\top\in\mathbb{R}^{n\times m} $, where the diagonal matrix $\Sigma^{\dagger}$
has the shape of the transpose of the diagonal matrix $\Sigma$, and is derived from $\Sigma$
by taking reciprocals of the non-zero (diagonal) elements of $\Sigma$ (i.e., the non-zero
singular values of $A$).

The following theorem gives a fundamental characterization of the M-P pseudoinverse.
\begin{theo}[\cite{Penrose}]
	For $A\in\mathbb{R}^{m \times n}$, the M-P pseudoinverse $A^{\dagger}$ is the unique $H\in\mathbb{R}^{n \times m}$ satisfying:
	\begin{align}
	& AHA = A \label{property1} \tag{P1}\\
	& HAH = H \label{property2} \tag{P2}\\
	& (AH)^{\top} = AH \label{property3} \tag{P3}\\
	& (HA)^{\top} = HA \label{property4} \tag{P4}
	\end{align}
\end{theo}

The first three M-P properties are particularly important for our purposes.

\medskip
\noindent {\bf \ref{property1}:} Following \cite{RohdeThesis}, we say that a \emph{generalized inverse} of $A$ is any $H$ satisfying  \ref{property1}.  Note that without \ref{property1}, even the all-zero matrix (which carries no information about $A$) satisfies the
other M-P properties. 
%The property \ref{property1} is particularly important in our context; without it,
%the all-zero matrix --- extremely sparse and carrying no information at all about $A$ --- would
%satisfy the other three properties.

\medskip
\noindent  {\bf \ref{property2}:} A generalized inverse is \emph{reflexive} if it satisfies \ref{property2}.  Theorem 3.14 in \cite{RohdeThesis} states that: ($i$) if $H$ is a generalized inverse of $A$, then $\mathrm{rank}(H)\ge\mathrm{rank}(A)$, and ($ii$) a generalized inverse $H$ of $A$ is reflexive if and only if $\mathrm{rank}(H)=\mathrm{rank}(A)$. Therefore, enforcing \ref{property2}
gives us the lowest possible rank of a generalized inverse of $A$.

\medskip
\noindent  {\bf \ref{property3}:} Following \cite{XuFampaLee},  we say that $H$ is \emph{ah-symmetric} if it satisfies \ref{property3}. That is, ah-symmetric means that $AH$ is symmetric.
If $H$ is an ah-symmetric generalized inverse (i.e., $H$ satisfies \ref{property1} and \ref{property3}),
then $\hat{x}:=Hb$ solves the least-squares problem $\min\{\|Ax-b\|_2:~x\in\mathbb{R}^n\}$.
So,  not  all of the M-P properties are required for a generalized inverse to solve this key problem. Especially important in our context, we have that by imposing property \ref{property3} on the 
generalized inverse when solving linear systems, we guarantee least-square solutions. If additionally,  \ref{property2} is imposed, we have the advantage of working with generalized inverses of minimum rank. Low rank is a desirable property for an ah-symmetric generalized inverse $H$,
in the context of the least-squares application, as it corresponds to 
a type of ``explainability'' for the associated linear model $\hat{x}:=Hb$. 

\medskip

Even if a given input matrix $A$ is sparse, its M-P pseudoinverse can be dense, leading to a high computational burden, especially in applications with  high-dimensional matrices.
  In \cite{FampaLee2018ORL,XuFampaLee,FLPX}, we propose and experiment with local-search procedures to construct 
 (``block-structured'')  reflexive generalized inverses having low 1-norm; effectively imposing minimum rank as a hard constraint.
 The use of the 1-norm as a surrogate for sparsity  is commonly applied across
 a wide variety of applications
 (see, for example, \cite{dokmanic1}) and is important in our context,
 so that we can construct  ah-symmetric generalized inverses using 
 tractable optimization problems. 
 Furthermore, in our context, use of the 1-norm has the very 
 important advantage of keeping the magnitude of the entries of the constructed generalized inverse under control
 --- this is extremely important for numerical stability in its applications.

 A 1-norm minimizing generalized inverse will not generally have minimum rank.
 A minimum-rank generalized inverse (i.e., a reflexive generalized inverse) will not generally
 have minimum 1-norm. So we seek generalized inverses that balance these two objectives. 
 In what follows, we investigate algorithmic options and outcomes associated with
 partially relaxing the reflexive property \ref{property2}, trading off 1-norm 
 against satisfaction of \ref{property2}, and in effect also trading off sparsity against rank.  
 
 %In \cite{FLPX}, numerical experiments with these procedures are presented. 
%  The purpose of the procedures is the construction of sparser matrices than the M-P pseudoinverse, 
% without losing some of its important properties. 
%In   \cite{FampaLee2018ORL,XuFampaLee},  (vector) 1-norm minimization is used to induce sparsity.  Therefore, at each iteration of the local-search procedures, the overall goal is to decrease the 1-norm of the constructed matrix $H$.\\

% Both sparsity and low rank are desirable features for ah-symmetric generalized inverses.  
% So, in this work, we are interested  in analyzing the trade-off between sparsity and low rank of these matrices, inducing sparsity  by (vector) 1-norm minimization. 

To obtain a $1$-norm minimizing ah-symmetric  generalized inverse $H$ of a given matrix $A$,  we  consider the optimization problem
\begin{equation}
\label{FPL}
\min\{\|H\|_1 \,:\, \ref{property1}, \ref{property3}\},
\end{equation}
where we write $\|H\|_1$ to mean $\|\mathrm{vec}(H)\|_1$, i.e.,  we use vector-norm notation on matrices.

To obtain  an ah-symmetric generalized inverse  with small $1$-norm, but also with guaranteed least rank, we should add  \ref{property2} to the constraints in (\ref{FPL}), resulting in
\begin{equation}
\label{FPL2}
\min\{\|H\|_1 \,:\, \ref{property1}, \ref{property2}, \ref{property3}\}.
\end{equation}

The solution of (\ref{FPL2}) is a $1$-norm minimizing ah-symmetric  reflexive generalized inverse of $A$. Clearly the price we pay for having the least rank condition on $H$, is a possible increase  in its $1$-norm, and we investigate through different approaches, how both measures, the $1$-norm and the rank of $H$, vary as we iteratively add  the individual constraints that model   \ref{property2} to problem (\ref{FPL}). 

We propose  algorithmic procedures that compute ah-symmetric generalized inverses of varied ranks and norms for a given $m\times n$ matrix $A$ of rank $r(A)<\min\{m,n\}$. The ah-symmetric generalized inverses constructed vary in a range that goes from matrices with minimum 1-norm to matrices with higher 1-norm and rank $r(A)$, i.e., the minimum possible rank. Through numerical experiments, we investigate if with the  procedures proposed, we can see a gradual  decrease in the rank as the norm increases. Intermediate solutions obtained by these algorithms can  realize a good trade-off between low 1-norm and low rank in applications of the generalized inverses. 

Concerning notation, in what follows, we use $J$ for an all-ones matrix, $r(A)$ for the rank of a matrix $A$, and  $\langle X, Y\rangle=\mathrm{trace}(X^\top Y):=\sum_{ij}x_{ij}y_{ij}$ for the matrix dot product.

%In what follows, for succinctness, we use vector-norm notation on matrices: we write $\|H\|_1$ to mean $\|\mathrm{vec}(H)\|_1$, and $\|H\|_{\max}$ to mean $\|\mathrm{vec}(H)\|_{\max}$ (in both cases, these are not the usual induced/operator matrix norms). We use $I$ for an identity matrix and $J$ for an all-ones matrix. Matrix dot product is indicated by $\langle X, Y\rangle=\mathrm{trace}(X^\top Y):=\sum_{ij}x_{ij}y_{ij}$. We use $A[S,T]$ for the submatrix of $A$ with row indices $S$ and column indices $T$; additionally, we use $A[S,:]$ ( resp., $A[:,T]$) for the submatrix of $A$ formed by the rows $S$ (resp., columns $T$).
%Finally, if $A$ is symmetric and $S=T$, we use $A[S]$ to represent the principal submatrix of $A$ with  row/column indices $S$.

We note that in order to formulate  (\ref{FPL2}) as a linear-programming (LP) problem, we linearize \ref{property2} using the following result.

\begin{prop}(see, for example, \cite[Proposition 4.3]{FFL2016})
	\label{prop1}
	If $H$ satisfies \ref{property1} and \ref{property3}, then $AH=AA^{\dagger}$.
\end{prop}
%\begin{proof}
%\begin{eqnarray*}
%AHA &=& AA^\dagger A \qquad \mbox{(by \ref{property1})}\\
%H'A'A &=& (A^\dagger )'A'A \qquad \mbox{(by \ref{property3})}\\
%A'AH &=& A'AA^\dagger \\
%(A^\dagger )'A'AH &=& (A^\dagger )'A'AA^\dagger \\
%AH&=&AA^\dagger ,
%\end{eqnarray*}
%the last equation following directly from a well-known property of $A^\dagger $.
%\end{proof}
Therefore, if  $H$ satisfies \ref{property1} and \ref{property3},   then 
\ref{property2} can be reformulated as the linear equation
\begin{equation}\label{p2lin} 
HAA^{\dagger}=H.
\end{equation}
Considering \eqref{p2lin}, we formulate  (\ref{FPL2}) as the linear program 
\[
\begin{array}{lll}
(P_{123}) \ z_{P_{123}}:=&\min  & \left\langle J , Z \right\rangle~,  \\
& \mbox{s.t.:}
&Z - H \geq 0~,\\
&&Z + H \geq 0~,\\
&&AHA=A~,\\
&&(AH)^{\top}=(AH)~,\\
&&HAA^{\dagger}=H~.
\end{array}
\]

We investigate  the trade-off  between   $1$-norm minimization and  low rank of $H$, by iteratively imposing \ref{property2} to problem (\ref{FPL}). We initially omit all constraints that model \ref{property2} in $P_{123}$, i.e., we solve the linear program 
\[
\begin{array}{lll}
(P_{13}) \ z_{P_{13}}:=&\min  & \left\langle J , Z \right\rangle~,  \\
& \mbox{s.t.:}
&Z - H \geq 0~,\\
&&Z + H \geq 0~,\\
&&AHA=A~,\\
&&(AH)^{\top}=(AH)~.
\end{array}
\]
Then, in different ways, we will gradually impose the ($mn$ individual) constraints  \eqref{p2lin}.

In \S\ref{sec:approaches}, we detailed describe our algorithmic approaches. In \S\ref{sec:numres}, we discuss our numerical results, and in \S\ref{sec:conclusions} we present our conclusions.

% For the purpose of comparison, we solve $P_{123}$ to obtain a minimum $1$-norm ah-symmetric reflexive generalized inverse. 
% For that, we use the solver Gurobi in two different ways. Considering the formulation of $P_{123}$, as it is, and considering the constraints that model \ref{property2} as \emph{lazy constraints}. %We note that, for this second approach, we have to add an integer variable to the problem, with no effect on its solution.  

\section{Algorithmic approaches}\label{sec:approaches}

In this section, we investigate different algorithmic approaches  to gradually impose the satisfaction of \ref{property2} on the solution of $P_{13}$. The purpose of our algorithms is to construct a set of diverse ah-symmetric generalized inverses, trading off low 1-norm against low rank. We aim for a sequence of intermediate ah-symmetric generalized inverses, constructed by our algorithms, trending toward lower rank at the expense of increasing 1-norm.  We propose five approaches and compare their ability to construct good 
solutions, where we obtain a  gradual decrease in the rank as the 1-norm iteratively increases.

The first four proposed algorithms use an interesting feature that is specific to our problem, of constructing  ah-symmetric generalized inverses: the equivalence between satisfaction of \ref{property2} and the least-rank condition. 
While \ref{property2} is a nonlinear equation in $H$ (hence beyond the realm of linear programming),
Proposition \ref{prop1} allows us to investigate the trade-off between 1-norm and rank by gradually imposing the linear equations \eqref{p2lin} on the linear program  $P_{13}$. 
We propose standard mathematical-programming approaches to gradually enforce these equations, namely, a cutting-plane method, an augmented Lagrangian method, and two penalty methods, where we penalize the 1-norm and the Frobenius norm of the violation,
given by the matrix $HAA^\dagger - H$, of the equations \eqref{p2lin}. The parameters adopted in these algorithms define both their initial solutions and how fast they converge to the ultimate solutions. When tuning these parameters, our objective is not to have a fast convergence.  On the contrary, we aim to construct a nice set of  ah-symmetric generalized inverses
%with varying 1-norms and ranks, 
approximating the Pareto frontier corresponding to our two minimization objectives,
1-norm and rank.  We design our algorithms to construct  a minimum 1-norm ah-symmetric generalized inverse of a given matrix $A$ at their first iteration, and we want them to slowly converge to a minimum 1-norm ah-symmetric \emph{reflexive} generalized inverse (which is necessarily a  minimum 1-norm ah-symmetric generalized inverse having minimum rank).
Along the way, by gradually imposing \ref{property2},  we aim to get a sequence of 
ah-symmetric generalized inverses of decreasing rank, with each having low 1-norm relative to  other
ah-symmetric generalized inverses of its rank.  

If not considering the particularity of our problem, given by the equivalence between \ref{property2} and the least-rank condition, a standard approach to balance objectives of having  low  1-norm and low rank when constructing a matrix, is to solve a convex optimization problem where the nuclear norm is employed
as a surrogate for the rank of the matrix, and then the objective is a weighted combination
of 1-norm and nuclear norm.  This  problem can be recast as a  semidefinite programming  (SDP) problem and has been widely investigated in the literature (see \cite{Parrilo,Recht,Fazel}, for example).
We also apply this general approach to iteratively construct  ah-symmetric generalized inverses.  However, compared to the four approaches discussed above, this method has the  disadvantage of requiring the solution of an SDP problem at each iteration, which generally does not scale well. Furthermore,  using  the nuclear norm
as a surrogate for rank, we lose a nice feature of our other approaches;
while the nuclear-norm approach will converge to an ah-symmetric  generalized inverse that is reflexive,
we lose the guarantee that it has minimum 1-norm among all such 
 ah-symmetric reflexive generalized inverses. %When applying this last approach, we adopt as the stopping criterion for the iterative algorithm proposed, the  condition  that the rank of $H$ is approximately equal to $r$. More specifically, we stop the algorithm when only $r$ singular values of $H$ are greater than $10^{-5}$. 

\subsection{Cutting-plane method for \ref{property2}}\label{sec:CP}
We  propose  a cutting-plane method, where  a linear program is solved at each iteration.
We initially solve problem $P_{13}$, not imposing \ref{property2}.  Then, we consider the $nm$ equations  \eqref{p2lin},  i.e., $H_{i\cdot} A (A^\dagger)_{\cdot j} = H_{ij}$,
lexically 
ordered by their indices $(j,i)$, with $i=1,\ldots,n$ and $j=1,\ldots,m$.
%corresponding to the elements $H_{ij}$ of the ah-symmetric generalized inverse constructed.  
At each iteration of the cutting-plane method, we add
the   first  $t$ equations that are violated by the current solution, to $P_{13}$.   We consider that equation $(j,i)$ is violated, if $|(HAA^\dagger - H)_{ij}| > 10^{-6}$. Aiming at a slow convergence of the algorithm, and therefore at a diverse set of constructed matrices $H$,  we limit  $t$  to 1\% of the total number of constraints.

\subsection{Augmented Lagrangian method: dualizing \ref{property2}}\label{sec:auglag}
Here, we initially consider a Lagrangian method, where we dualize the constraints  \eqref{p2lin}.   We solve the Lagrangian-dual problem
\begin{equation}\label{lag-dual}
\begin{array}{ll}
\max_{\Lambda}\min_{H,Z}  & \left\langle J , Z \right\rangle + \left\langle \Lambda , HAA^{\dagger}-H \right\rangle~,  \\
 \mbox{s.t.:}
&Z - H \geq 0~,\\
&Z + H \geq 0~,\\
&AHA=A~,\\
&(AH)^{\top}=(AH)~,%\\
%&&HAA^{\dagger}-H=K~,
\end{array}
\end{equation}
with a subgradient optimization algorithm. At each iteration, we solve the inner problem  in \eqref{lag-dual} fixing the Lagrangian  multiplier $\Lambda\in\mathbb{R}^{n\times m}$ at a value $\Lambda_k$ and obtain its solution $(H^k,Z^k)$. We calculate  the subgradient $G_k:=H^kAA^{\dagger}-H^k$ and  update the Lagrangian multiplier according to  $ \Lambda_{k+1}:=\Lambda_k + \gamma_k G_k$. The step size $\gamma_k$ is $(z_{P_{123}} - z_k)/\|H^kAA^{\dagger}-H^k\|_F$, where $z_k$ is the current value of the Lagrangian function, i.e., of the objective function  in \eqref{lag-dual}.  
% We  initialize the algorithm with $\Lambda=\mathbf{0}$. 
% The  algorithm stops if all the absolute values of the elements of $K$ are less than a tolerance $\epsilon$, or if  the number of iterations reaches  $1000$.

We note that although we have a theoretical guarantee of convergence of the subgradient algorithm to the optimal value of $P_{123}$ when applying the Lagrangian method, there is no guarantee that a feasible solution for $P_{123}$ is obtained. We actually noticed in preliminary numerical experiments that the rank of $H$ did not converge to the rank of $A$, i.e.,  \ref{property2} was not satisfied by the solution obtained when the algorithm converged. To overcome this drawback, we have considered the augmented Lagrangian method,  adding the convex quadratic term  $\frac{\mu_k}{2} \|HAA^{\dagger}-H\|_F$ to the objective function of the Lagrangian subproblem. At  iteration $k$ of the algorithm, we solve the convex quadratic  problem 
\[
\begin{array}{ll}
\min_{H,Z}  & \left\langle J , Z \right\rangle + \left\langle \Lambda_k , HAA^{\dagger}-H \right\rangle  + \frac{\mu_k}{2} \|HAA^{\dagger}-H\|_F~,  \\
 \mbox{s.t.:}
&Z - H \geq 0~,\\
&Z + H \geq 0~,\\
&AHA=A~,\\
&(AH)^{\top}=(AH)~,%\\
%&HAA^{\dagger}-H=K~,
\end{array}
\]
for  fixed $\mu_k$ and $\Lambda_k$, obtaining its optimal solution $(H_k,Z_k)$. After solving the problem, we update the penalty parameter  according to: $\mu_{k+1}=1.30 \mu_{k}$ and then we update the dual variable according to $\Lambda_{k+1}:=\Lambda_{i} + \mu_{k+1} (H_kAA^{\dagger}-H_k)$.
We initialize the algorithm with $\Lambda_0:=0$ and $\mu_0:=0.1$.

\subsection{Penalty method: Frobenius norm of \ref{property2} violation}%penalizing \texorpdfstring{$\|HAA^{\dagger}-H\|_F$}{HAA-HF}  }  
Now we consider penalty methods, penalizing the Frobenius norm of $HAA^{\dagger}-H$. At iteration $k$ of the algorithm,  we solve the convex quadratic problem
\[
\begin{array}{ll}
\min_{H,Z}  & \left\langle J , Z \right\rangle + \mu_k \|HAA^{\dagger}-H\|_F~,  \\
 \mbox{s.t.:}
&Z - H \geq 0~,\\
&Z + H \geq 0~,\\
&AHA=A~,\\
&(AH)^{\top}=(AH)~,%\\
%&&HAA^{\dagger}-H=K~,
\end{array}
\]
for  fixed $\mu_k$. After solving the problem, we update the penalty parameter  according to: $\mu_{k+1}:=1.30 \mu_{k}$.
We   initialize the algorithm with  $\mu_0:=0.1$.

\subsection{Penalty method: 1-norm of \ref{property2} violation} %{Penalty method, penalizing \texorpdfstring{$\|HAA^{\dagger}-H\|_1$}{HAA-HF} } 
We also experimented penalizing the $1$-norm of  $HAA^{\dagger}-H$. Besides the modification on the norm with respect to the  penalty method described in the previous subsection, in this algorithm, we have tuned the parameters to $\mu_0:=0.01$ and $\mu_{k+1}:=1.15 \mu_{k}$. We note that the use of the $1$-norm instead of the Frobenius norm leads to the solution of a linear program at each iteration, instead of a convex quadratic program. 

\subsection{Nuclear-norm method}\label{sec:penrank}
An standard approach to obtain a balance between the two goals of constructing a sparse matrix and a low-rank matrix,  is to minimize a weighted sum of its 0-norm and  rank. In our context this would lead to the following problem 
\begin{equation}\label{0-rank}
\min\{\|H\|_0 + \mu \cdot  r(H) \,:\, \ref{property1},\ref{property3}\},
\end{equation}
where $\mu$ corresponds to the weight given to the rank of $H$ ($r(H)$) in the minimization. 
To address a convex optimization problem, it is usual then to employ the 1-norm as a surrogate for the 0-norm, as we have done with the previous approaches proposed,  and the nuclear norm as a surrogate for the rank as well. Applying this method to our problem, we solve at iteration $k$ of our last proposed algorithm, the following  problem
\begin{equation}
\label{rankSDP}
\min\{\|H\|_1 + \mu_k \|H\|_* \,:\, \ref{property1},\ref{property3}\},
\end{equation}
where $\mu_k$ is a fixed real penalty parameter and $\|H\|_*:=\sum_k\sigma_k(H)$ (the nuclear norm), where $\sigma(H) := (\sigma_1(H), . . . , \sigma_n(H))$
is the vector of singular values of $H$. It is convenient and usual to assume that the singular
values are ordered so that $\sigma_1(H) \geq \ldots \geq \sigma_n(H)$. After solving the problem, we update the penalty parameter  according to $\mu_{k+1}:=1.20 \mu_{k}$.
We initialize the algorithm with  $\mu_0:=0.05$. 

Problem (\ref{rankSDP}) can be recast as an SDP problem. Using
the facts that the nuclear norm is dual to the spectral norm and that the spectral norm admits
a simple semidefinite characterization, minimizing $\|H\|_*$ can be formulated as
\[
\min_{W_1,W_2}\left\{\frac{1}{2}(\mbox{trace}(W_1)+\mbox{trace}(W_2)) \,:\,\left(\begin{array}{cc}W_1&H\\H^T&W_2\end{array}\right)\succeq 0\right\}.
\]
Therefore, at each iteration $k$ of our algorithm, we solve the SDP problem
\begin{equation}\label{sdp_prob}
\begin{array}{l}
\displaystyle \min_{H,Z,W_1,W_2}\left\langle J , Z \right\rangle +\mu_k\left( \frac{1}{2}(\mbox{trace}(W_1)+\mbox{trace}(W_2))\right)\\
\mbox{subject to:}\\
\left(\begin{array}{cc}W_1&H\\H^T&W_2\end{array}\right)\succeq 0,\; -Z\leq H\leq Z, \; AHA=A,\; (AH)^\top=(AH),
\end{array}
\end{equation}
where the matrix variables $W_1$ and $W_2$ are of order $n\times n$ and $m\times m$, respectively. 

\section{Numerical results}\label{sec:numres}

To analyze the results of the procedures proposed, we randomly generated 20 dense square matrices of size $50\times 50$ and  rank $r=25$.
We used the Matlab function \emph{sprand}, which  generates a random  $m\times n$ dimensional matrix $A$ with singular values
 given by a nonnegative input vector $rc$. 
 We  selected the $r$ nonzeros of $rc$
 as the decreasing vector $M\times(\rho^{1},\rho^{2},\ldots,\rho^{r})$, where $M=2$, and $\rho=(1/M)^{(2/(r+1))}$.
We implemented our  algorithms in Python (Spyder 3.3.6), and ran the experiments on a
16-core machine (running Windows Server 2016 Standard):
two Intel Xeon CPU E5-2667 v4 processors
running at 3.20GHz, with 8 cores each, and 128 GB of memory. We solve all the optimization problems involved in our experiments  with Gurobi v.9, except problem \eqref{sdp_prob}, which we solve with Mosek.

In our numerical analysis, when we refer to the rank and $0$-norm of a generalized inverse computed by our algorithms, we mean, respectively,  the number of singular values  greater than $10^{-5}$   and the number of elements of the matrix with absolute values  greater than $10^{-6}$. When we refer to the number of constraints  in \eqref{p2lin} that are satisfied by a given matrix $H$, we  mean the number of elements of the matrix $HAA^\dagger - H$ with absolute value less than $10^{-6}$.  We denote  optimal solutions of $P_{13}$ and $P_{123}$, respectively by $H_{13}$ and $H_{123}$.  Most of the results presented in this section are average results over our 20 test-instances.

%  We present results for all algorithms described in the previous section except for the Lagrangian method, which was disregarded because, although its optimal solution value converges to $\|H_{123}\|_1$, the rank of $H$ does not converge to the rank of $A$. Moreover, the 1-norm of $H$ oscillates very much during the execution of the algorithm. For this reason, we have replaced this algorithm with the Augmented Lagrangian algorithm, described in Subsection \ref{sec:auglag}. The parameters used in the algorithms are all given in the previous section. They were tuned  with the goal of having all the algorithms starting with approximately the same solution ($H_{13}$), and converging with approximately the same number of iterations, except penalized rank, which stops with less iterations. 

In Table \ref{table1}, we analyze the solutions which the methods proposed in the previous section converge to, and compare them to the M-P pseudoinverse $A^\dagger$, to  $H_{13}$, and to  $H_{123}$. The stopping criterion for all of our algorithms is the same. They stop when for the constructed matrix $H$, we have $|(HAA^\dagger - H)_{ij}|\leq 10^{-6}$, for  $i=1,\ldots,n$ and $j=1,\ldots,m$. We present  in the table,   the average number of iterations executed by the algorithms,  
and the  average 1-norm, 0-norm and rank of the solutions obtained.  
We note that all constraints  in \eqref{p2lin} are satisfied by the solutions of the 
methods proposed, and therefore they all have the minimum rank $r=25$. 

\begin{table}[!ht]
\centering
\begin{tabular}{l|rrrc} 
\hline
                & \multicolumn{1}{c}{It}  & 
                %\multicolumn{1}{c}{Consts.} & %\multicolumn{1}{c}{P2}&
                %\multicolumn{1}{c}{Time} &  
                \multicolumn{1}{c}{$\|H\|_1$} &  \multicolumn{1}{c}{$\|H\|_0$} &\multicolumn{1}{c}{$r(H)$} \\ 
                \hline
$A^\dagger$      & \multicolumn{1}{c}{-}                % & \multicolumn{1}{c}{-}              %    & 100.0      
%& 0.005
& 157.2                     & 2384.0    &25     \\
 $H_{13}$         & \multicolumn{1}{c}{-}                 %& 2500.0            % & 100.0       
 %& 168.0       
 &           131.7           &    1218.6  &  43.8 \\
 $H_{123}$         & \multicolumn{1}{c}{-}                 %& 2500.0            % & 100.0       
 %& 168.0       
 & 137.9                     & 1581.6     &25    \\
\hline

Cutting-plane  & 45.7              %& 1110.7         %    & 100.0       %& 557.8    
& 137.5                 & 1596.6   &25      \\
Aug. Lagrangian  & 46.4             %& \multicolumn{1}{c}{-}                 % & 100.0       %& 511.4 
& 137.5                 & 1597.0    &25    \\
Pen. 1-norm      & 33.8              %& \multicolumn{1}{c}{-}              %    & 100.0      % & 2668.8    
& 137.5               & 1596.7    &25    \\
Pen. Frobenius & 60.8             % & \multicolumn{1}{c}{-}                 % & 100.0      % & 624.7    
& 137.5                  & 1597.4     &25    \\
Nuclear-Norm      & 23.4             % & \multicolumn{1}{c}{-}               %   & 100.0       % & 17534.9  
& 142.4                & 1897.7   &25 \\ \hline     
\end{tabular}
	\caption{Comparison of solutions of the methods proposed with $A^\dagger$, $H_{13}$, and  $H_{123}$ \label{table1}}
\end{table}

From the results in Table \ref{table1}, we  note that the M-P pseudoinverse has 
significantly greater 1-norm and density (0-norm) than  $H_{123}$, which indicates a significant advantage of using our 
minimum 1-norm ah-symmetric reflexive generalized inverse in numerical applications. 
In fact, the M-P pseudoinverse has about 95\% nonzero elements, on average, while $H_{123}$ has only about 63\%.  Comparing $H_{13}$ to $H_{123}$, we see that when imposing the reflexive property \ref{property2} on the ah-symmetric generalized inverse, we obtain  an average  decrease of 43\% in the rank of the matrix at the expense of increasing the 1-norm and the 0-norm by
approximately 4.5\% and 30\%, respectively. In summary, if we insist on minimum rank then we would use $H_{123}$.
But we can gain alternatively use $H_{13}$, gaining considerable decreases in 1-norm and density,
but with much greater rank. The space between these extreme alternatives is what we aim to explore.
We also observe from the results  that, except for the Nuclear-norm method,  the solutions $H$ obtained by all the methods have, on average, the same  1-norms and only slightly different 0-norms. The   1-norms are slightly less than   the  1-norm  of $H_{123}$,  because equations \eqref{p2lin} are only satisfied up to a tolerance. The 0-norms are only slightly greater than the 0-norm of $H_{123}$.  
The Nuclear-norm method gains full rank only at the expense of barely weighting the 1-norm. 

In the next experiment, we analyze intermediate solutions obtained by the methods  with the purpose of seeing their ability to construct solutions  with different 1-norms and ranks as \ref{property2} is gradually imposed. Our goal now, is to analyze the trade-off between the 1-norms and 0-norms of the intermediate solutions  and their ranks.  
 
In Table \ref{table2}, we present statistics for the first iteration of each  method in which the rank of the computed generalized inverse $H$  has sufficiently decreased, more specifically, for the first iteration where  the rank of $H$ satisfies
\begin{equation}\label{alpha}
r(H)\leq (1-\alpha) \cdot r(H_{13}) + \alpha \cdot r(H_{123}),\end{equation}
for $\alpha=0.00,0.25,0.50,0.75,1.00$. We tabulate the average percentage
%relative 
increase in the 1-norm and  0-norm of $H$ compared to the initial matrix $H_{13}$, i.e., 
\[
\widehat{\|H\|}_i:=\frac{\|H\|_i - \|H_{13}\|_i}{\|H_{13}\|_i}  \times 100 ,
\]
for $i=0,1$. For each $\alpha$, we also present $r(H)$, given by \eqref{alpha}, the average percentage of constraints in \eqref{p2lin} that are satisfied by the intermediate solution $H$ (`P2') and the average iteration in which it is computed (`It').  

\begin{table}[!ht]
\centering{
	\begin{tabular}{c|c|rrrr|rrrr}
		\hline
&&\multicolumn{4}{c|}{Cutting-plane}&\multicolumn{4}{c}{Aug. Lagrangian}\\
$\alpha$&$r(H)$&$\widehat{\|H\|}_1$&$\widehat{\|H\|}_0$&P2&It&$\widehat{\|H\|}_1$&$\widehat{\|H\|}_0$&P2&It\\ \hline
0.00&43.8&0.0&0.0&14.8&0&0.0&0.0&14.9&0\\
0.25&39.1&2.5&15.1&36.2&18&4.3&29.0&50.2&31\\
0.50&34.4&3.4&21.5&46.4&27&4.3&29.4&69.1&33\\
0.75&29.7&4.1&28.1&59.6&36&4.3&29.8&87.5&36\\
1.00&25.0&4.3&31.0&100.0&51&4.3&30.2&100.0&51\\
	\end{tabular}
		\begin{tabular}{c|c|rrrr|rrrr|rrrr}
		\hline
&&\multicolumn{4}{c|}{Pen. 1-norm}&\multicolumn{4}{c|}{Pen. Frobenius}&\multicolumn{4}{c}{Nuclear-Norm}\\
$\alpha$&$r(H)$&$\widehat{\|H\|}_1$&$\widehat{\|H\|}_0$&P2&It&$\widehat{\|H\|}_1$&$\widehat{\|H\|}_0$&P2&It&$\widehat{\|H\|}_1$&$\widehat{\|H\|}_0$&P2&It\\ \hline
0.00&43.8&0.0&0.0&15.4&0&0.0&0.0&14.9&0&0.0&0.0&14.7&0\\
0.25&39.1&4.3&29.6&68.4&29&4.3&30.4&50.5&50&3.1&24.6&12.7&17\\
0.50&34.4&4.3&29.9&84.5&30&4.3&30.3&71.3&53&4.8&34.7&12.1&19\\
0.75&29.7&4.3&30.0&92.4&31&4.3&30.3&85.7&55&6.5&45.3&11.5&21\\
1.00&25.0&4.3&30.1&100.0&35&4.3&30.2&100.0&61&8.1&53.7&100.0&24\\
\hline
	\end{tabular}
	\caption{Trade-off between norms and rank for ah-symmetric generalized inverses \label{table2}}
}
\end{table}
From the results in Table \ref{table2}, we observe that 
\begin{itemize}
\item[$\bullet$] the Cutting-plane method presents the best trade-off between  norms and rank of $H$.  Approximately 35\%  of the iterations are executed before the algorithm finds the first matrix $H$ with rank satisfying \eqref{alpha} for $\alpha=0.25$, but matrices with different norms are computed for each $\alpha=0.0,0.25,0.50,0.75$. Both 1-norm and 0-norm increase as $\alpha$ gradually increase   in this interval; 

\item[$\bullet$] the methods Augmented Lagrangian, Penalized 1-norm, and Penalized Frobenius have similar behavior. 
For the values of $\alpha$ considered in Table \ref{table2}, they only start to decrease the rank of $H$ when the 1-norm reaches its maximum value. Therefore, for the values of $\alpha$ considered in this experiment, the solutions generated by these methods do not show the desired trade-off between 1-norm and rank;

\item[$\bullet$] on average, the 1-norm and the 0-norm of the final solutions of all methods, except Nuclear-norm, are, respectively,   4\% and 30\% greater than the norms of the initial matrix $H_{13}$. We can see that the 0-norm increases as the 1-norm increases, demonstrating that the 1-norm works as a good surrogate for the 0-norm in this experiment;

\item[$\bullet$] all methods, except Nuclear-norm, converge  to matrices $H$ with similar norms that satisfy all  the constraints in \eqref{p2lin};

\item[$\bullet$] Nuclear-norm converges to worse solutions than the other methods proposed. The solutions have greater norms on average. Unlike we do in the other methods, we do not enforce the constraints \eqref{p2lin} in this one. Therefore,  the number of constraints in \eqref{p2lin} that are satisfied, is very small during the execution of the method, it increases very fast only in the last iterations, when the rank of $H$ reaches its lowest value.  We see that the algorithm converges  to a solution $H$ of least rank, but with $1$-norm greater than  $\|H_{123}\|_1$, which is due to the fact that \eqref{sdp_prob} is not a relaxation for $P_{123}$, in contrast to the problems solved by  the other methods.   
Although we can observe a gradual increase in the norms for this method, because of the higher norms, it is not a good option to generate intermediate ah-symmetric generalized inverses. 
\end{itemize}

In conclusion, from the results in Table \ref{table2}, we see that among the first four methods, which  converge to matrices of similar norms, Cutting-plane is the only one that shows a gradual increase in the norms for the values of $\alpha$ selected. The others only start to decrease the rank of $H$ after its 1-norm  reaches its greatest value. Nuclear-norm also shows a gradual increase in the norm, however, it converges to solutions of higher norms (on average). Therefore,  our experiment points to Cutting-plane as the most suitable among the proposed methods to trade-off 1-norm against rank for ah-symmetric generalized inverses.

\begin{figure}[!ht]
	\centering
	\includegraphics[width=0.95\textwidth]{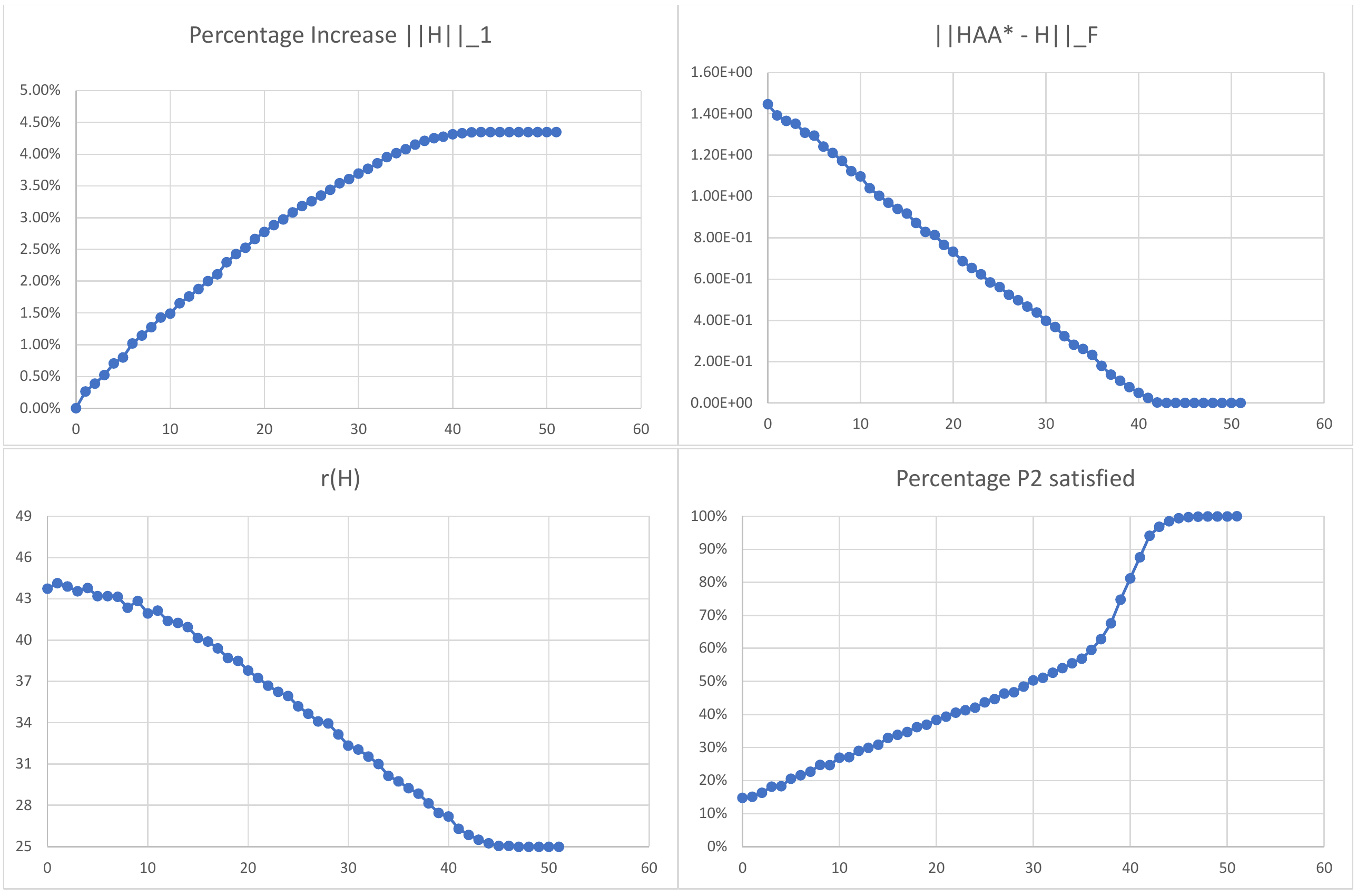}
	\caption{Cutting-plane method}\label{figure1a}
\end{figure}

\begin{figure}[!ht]
	\centering
    \includegraphics[width=0.95\textwidth]{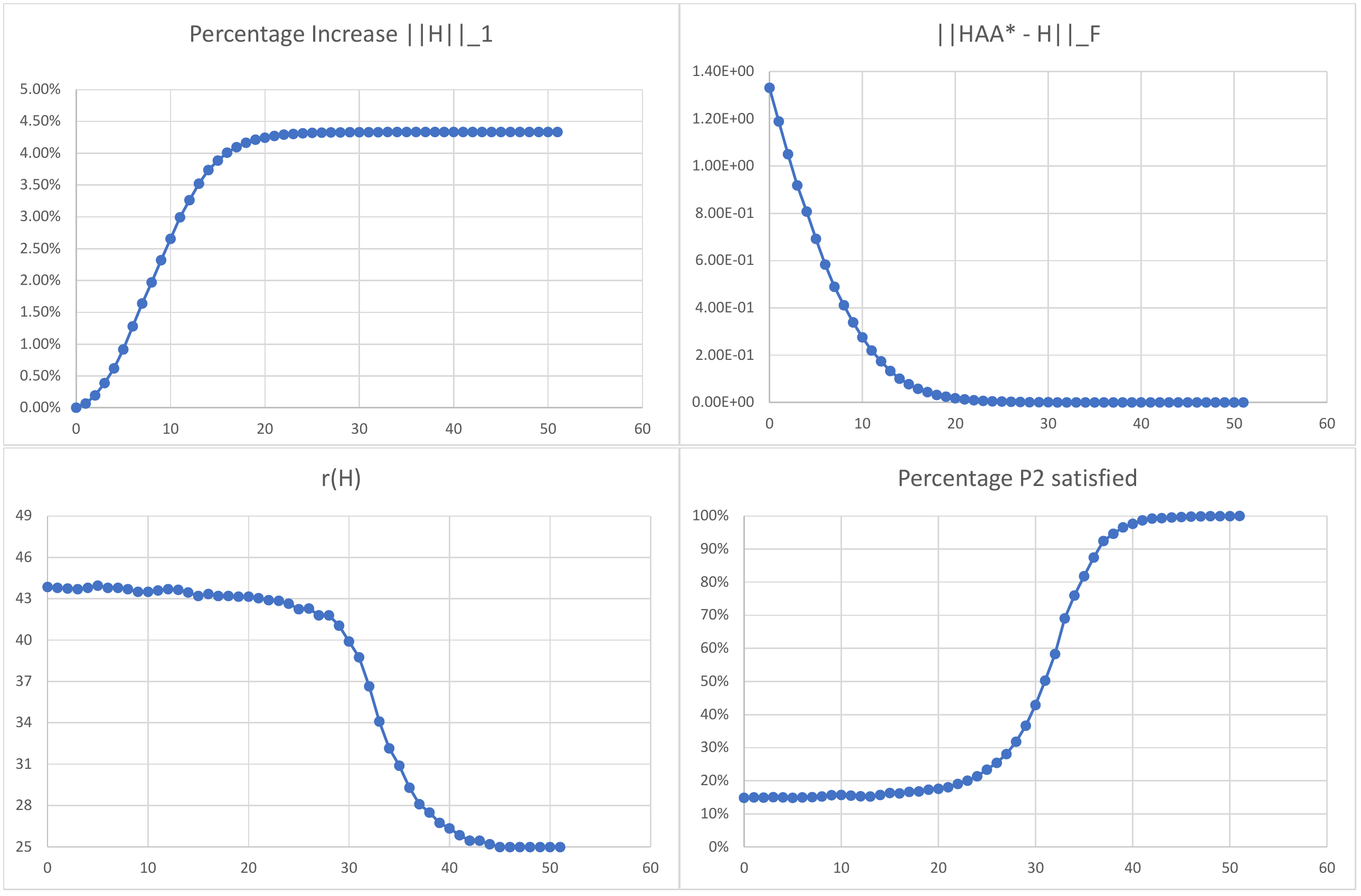}
	\caption{Augmented Lagrangian method}\label{figure1b}%-  $(50,50,25,1.00)$}
\end{figure}

\begin{figure}[!ht]
	\centering
	\includegraphics[width=0.95\textwidth]{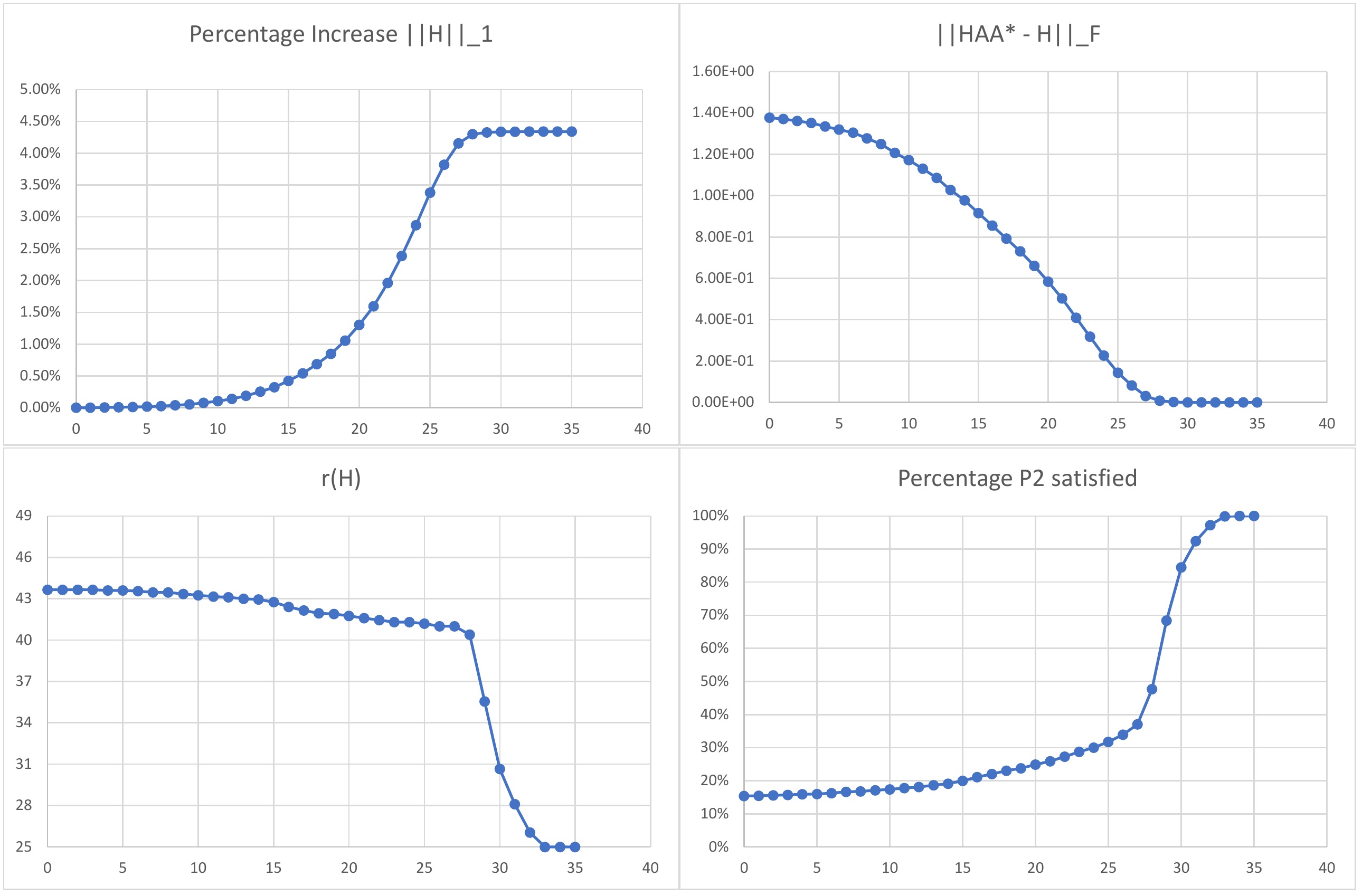}
	\caption{Penalized 1-norm method }\label{figure2a}
\end{figure}

\begin{figure}[!ht]
	\centering
	\includegraphics[width=0.95\textwidth]{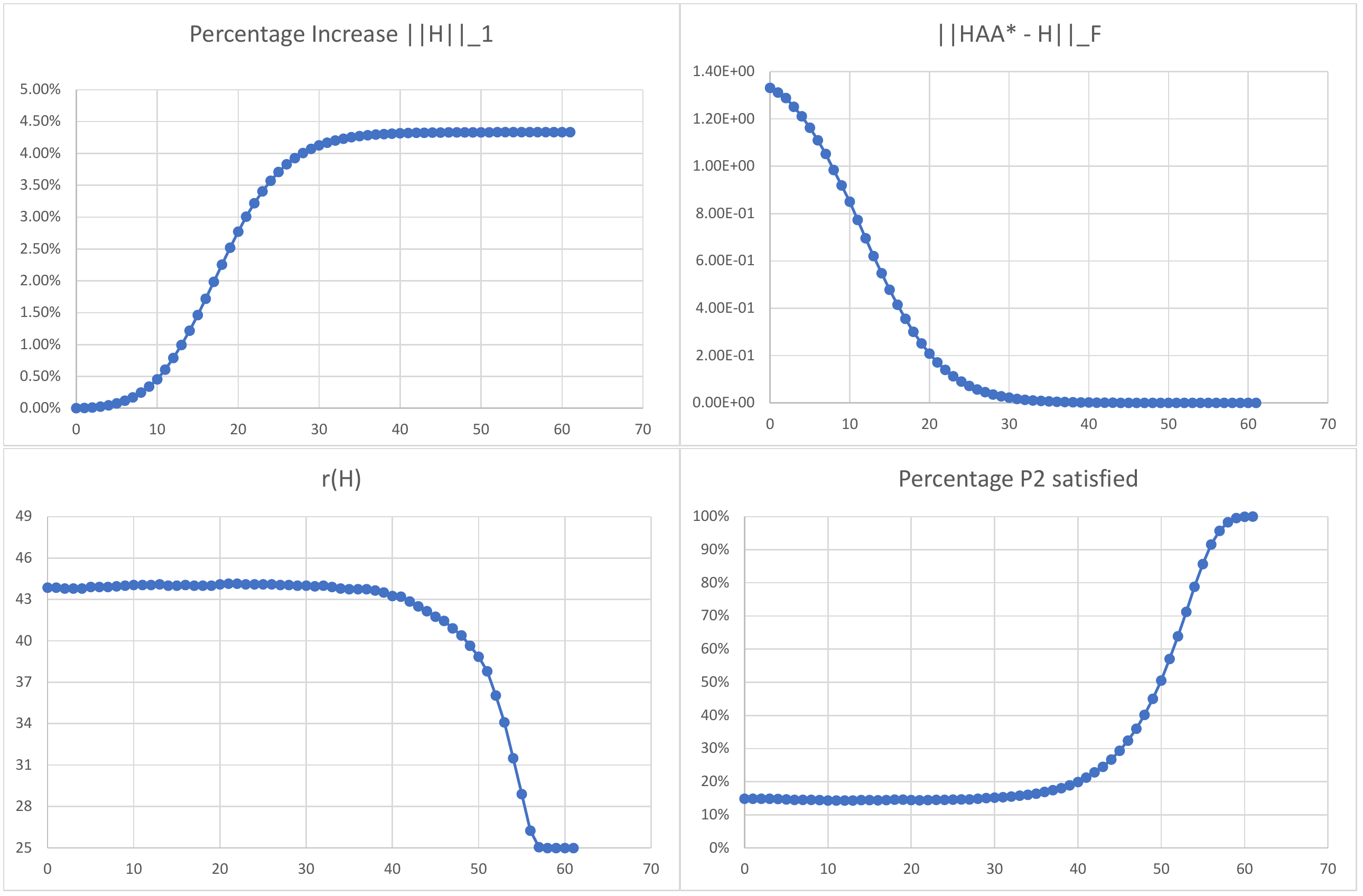}
	\caption{Penalized Frobenius method}\label{figure2b}
\end{figure}

\begin{figure}[!ht]
	\centering
    \includegraphics[width=0.95\textwidth]{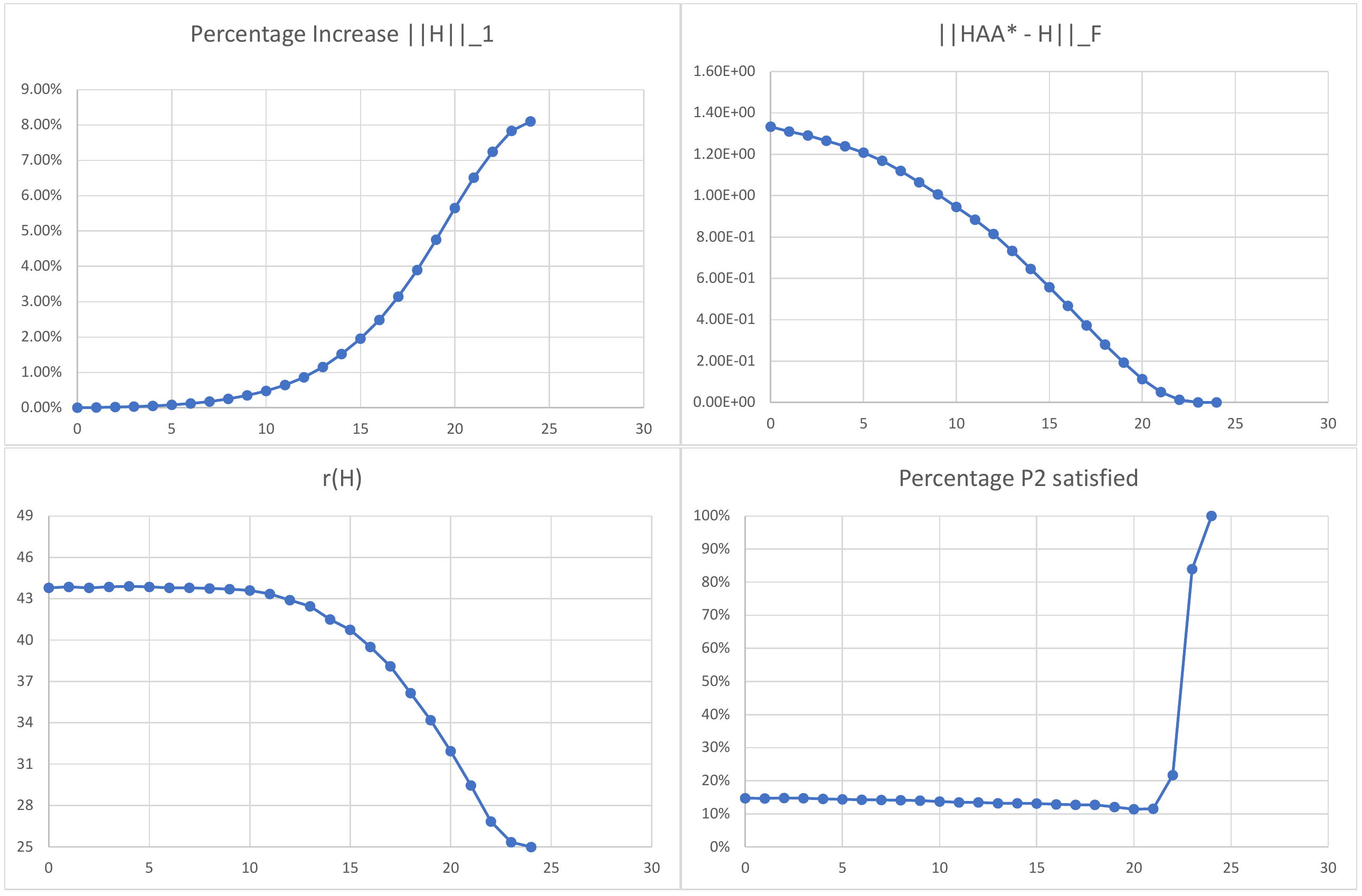}
	\caption{ Nuclear-norm method}\label{figure2c}
\end{figure}

In Figures \ref{figure1a} to \ref{figure2c}, we present plots to better show  the behavior of the  methods proposed as they iterate. 
The plots depict  average values for $\widehat{\|H\|}_1$,  
%$(z_{P_{123}}- \|H\|_1)/z_{P_{123}}$, 
$\|HAA^\dagger - H\|_F$, $r(H)$, and  the percentage of the constraints in \eqref{p2lin} that are satisfied. 

Comparing the plots for the Cutting-plane method, in Figure \ref{figure1a}, with the plots for Augmented Lagrangian, Penalized  1-norm and Penalized Frobenius, in  Figures \ref{figure1b}, \ref{figure2a}, \ref{figure2b}, we can more clearly observe that during the execution of Cutting-plane the rank starts to decrease when the 1-norm is still increasing and the percentage of constraints in \eqref{p2lin} that are satisfied also increases faster from the beginning. The behavior of the solutions of the three other methods is very similar.
% Although the 1-norm increases faster for Augmented Lagrangian,  even for the other two, the rank does not decrease significantly while  the 1-norm is still increasing. 
Nevertheless,  Penalized 1-norm shows a slightly better trade-off than the other two. We can observe a small decrease in the rank before the 1-norm achieves its optimal value. In Figure \ref{figure2c}, we see that Nuclear-norm converges in fewer iterations than the other methods, but in terms of computational effort, it is the most expensive method, owing to the need to solve SDP problems. 
As guaranteed by the theory, all of the constraints \eqref{p2lin} are satisfied when the rank is minimum, but only a few of them  are satisfied until the last iterations. Despite this fact, the curves show the desired  decrease in the rank as the 1-norm increases. However, as observed in Table \ref{table2} the 1-norm increases significantly more for this method than for the others. We have an increase of 8.1\% in the 1-norm, while for the others it is only 4.3\%. 

% first 
% O metodo Penalized Rank, e um método bem demorado, além disso, foi necessário utilizar um fator multiplicativo (stepsize) muito baixo, pois o problema se torna inviavel com um stepsize um pouco maior. Comecamos com stepsize pequeno para inicialmente satisfazer P13.

% Penalized Frobenius e Augmented Lagrangian sao muitos parecidos, isso pq o Lagrangian nao altera muito bem a solucao, a maior diferenca e no termo quadratico. Porem, o Augmented Lagrangian apresenta uma estabilidade numerica maior, isso deve ser pelo Lagrangian, o Augmented Lagrangian apresenta valores singulares melhor no final

\FloatBarrier

\FloatBarrier

As Cutting-plane presented the best trade-off between the rank and the 1-norm of the solutions in the experiments discussed above, we further investigated if other ways of selecting the violated constraints to add to the problem at each iteration could improve even more the results. The best results were obtained when the violated constraints were selected randomly. So, in the following, we include in our analysis, a variation of Cutting-plane called Cutting-plane Random. Recall that for the original Cutting-plane,  
we consider the $nm$ equations  \eqref{p2lin},  i.e., $H_{i\cdot} A (A^\dagger)_{\cdot j} = H_{ij}$,
lexically ordered by their indices $(j,i)$, with $i=1,\ldots,n$ and $j=1,\ldots,m$.
% we consider the $nm$ equations  \eqref{p2lin},  lexically 
% ordered by their indices $(j,i)$, with $j=1,\ldots,m$ and $i=1,\ldots,n$, corresponding to the elements $H_{ij}$ of the ah-symmetric generalized inverse constructed. 
At each iteration of Cutting-plane, we add
the   first  $t$ equations that are violated by the current solution, to $P_{13}$. Alternatively, in Cutting-plane Random, we randomly reorder the indices $j$ or the indices $i$ before executing it. We execute the method twenty times, in the first ten, we randomly reorder $j$ and in the last ten,  besides randomly reordering $i$, we search considering the lexical order of $(i,j)$. 
The solution of minimum 1-norm, for each rank found in the solutions is then presented. 
We will compare in the following Cutting-plane (CP), Cutting-plane Random (CP.Rand), Penalized 1-norm, and  Nuclear-norm. The two other methods, Augmented Lagrangian and Penalized Frobenius, are not considered in the following because their behavior is very similar to Penalized 1-norm, only a bit worse.

\begin{table}[!ht]
\centering{
	\begin{tabular}{c|cccc|cccc}
		\hline
		&\multicolumn{4}{c|}{$\|H\|_1$}&\multicolumn{4}{c}{$\|H\|_0$}\\
		$r(H)$&Nuclear-norm&Pen.1-norm&CP& CP.Rand&Nuclear-norm&Pen.1-norm&CP& CP.Rand\\
		\hline
48 & (1, 0)  & (1, 0)  & (1, 1)  & (1, 1)   & (1, 0)  & (1, 0)   & (1, 1)   & (1, 1)   \\
47 & (3, 0)  & (3, 0)  & (3, 2)  & (3, 3)   & (3, 0)  & (3, 0)   & (3, 2)   & (3, 3)   \\
46 & (5, 1)  & (5, 0)  & (5, 2)  & (5, 4)   & (5, 0)  & (5, 0)   & (5, 3)   & (5, 4)   \\
45 & (7, 0)  & (5, 1)  & (8, 3)  & (8, 7)   & (7, 0)  & (5, 0)   & (8, 4)   & (8, 7)   \\
44 & (11, 0) & (9, 2)  & (11, 3) & (11, 9)  & (11, 0) & (9, 0)   & (11, 5)  & (11, 10) \\
43 & (13, 0) & (8, 2)  & (15, 6) & (15, 12) & (13, 0) & (8, 0)   & (15, 8)  & (15, 13) \\
42 & (14, 1) & (11, 1) & (16, 3) & (16, 12) & (14, 0) & (11, 0)  & (16, 7)  & (16, 12) \\
41 & (15, 0) & (12, 2) & (19, 5) & (19, 15) & (15, 0) & (12, 0)  & (19, 7)  & (19, 15) \\
40 & (13, 1) & (9, 3)  & (19, 3) & (19, 12) & (13, 1) & (9, 1)   & (19, 8)  & (19, 11) \\
39 & (7, 0)  & (8, 0)  & (19, 5) & (19, 14) & (7, 0)  & (8, 0)   & (19, 8)  & (19, 12) \\
38 & (8, 0)  & (3, 0)  & (18, 2) & (19, 17) & (8, 0)  & (3, 0)   & (18, 4)  & (19, 15) \\
37 & (13, 0) & (4, 0)  & (20, 5) & (20, 16) & (13, 0) & (4, 0)   & (20, 9)  & (20, 13) \\
36 & (8, 0)  & (5, 0)  & (20, 3) & (20, 17) & (8, 0)  & (5, 0)   & (20, 7)  & (20, 13) \\
35 & (10, 0) & (5, 0)  & (20, 1) & (20, 19) & (10, 0) & (5, 0)   & (20, 7)  & (20, 13) \\
34 & (11, 0) & (4, 1)  & (20, 1) & (20, 18) & (11, 0) & (4, 0)   & (20, 7)  & (20, 13) \\
33 & (9, 0)  & (5, 1)  & (20, 1) & (20, 18) & (9, 0)  & (5, 0)   & (20, 8)  & (20, 13) \\
32 & (8, 0)  & (6, 1)  & (20, 1) & (20, 18) & (8, 0)  & (6, 1)   & (20, 8)  & (20, 12) \\
31 & (9, 0)  & (5, 0)  & (20, 1) & (20, 19) & (9, 0)  & (5, 0)   & (20, 8)  & (20, 12) \\
30 & (7, 0)  & (5, 0)  & (20, 1) & (20, 19) & (7, 0)  & (5, 0)   & (20, 5)  & (20, 16) \\
29 & (9, 0)  & (8, 0)  & (19, 4) & (20, 16) & (9, 0)  & (8, 0)   & (19, 11) & (20, 10) \\
28 & (4, 0)  & (12, 0) & (20, 3) & (20, 17) & (4, 0)  & (12, 0)  & (20, 5)  & (20, 17) \\
27 & (12, 0) & (6, 0)  & (20, 2) & (20, 18) & (12, 0) & (6, 1)   & (20, 6)  & (20, 14) \\
26 & (12, 0) & (12, 0) & (20, 3) & (20, 18) & (12, 0) & (12, 0)  & (20, 6)  & (20, 15) \\
25 & (20, 0) & (20, 2) & (20, 6) & (20, 18) & (20, 0) & (20, 19) & (20, 17) & (20, 14) \\
		\hline
		\end{tabular}
			\caption{Number of instances where the methods find (at least one solution, the solution of minimum norm) \label{table3}}
		}
		\end{table}

In Table \ref{table3}, each row corresponds to the rank of at least one solution generated by the four methods when applied to our twenty instances. For each method, we show results concerning the 1-norm and the 0-norm of the solutions obtained. Each ordered pair presented in the table, contains the number of instances where the methods find at least one solution for the corresponding rank, and the number of instances where the method finds the least 1-norm (up to a tolerance of $10^{-4}$), for that rank, among all of the algorithms.  A method for which the first coordinate is large indicates the ability of that method to generate
 a diverse set of solutions (i.e., with a variety of ranks), when trading-off rank against 1-norm. 
 The second coordinate indicates when a method obtain the best solution among all for that rank, hence 
giving a nondominated solution, in the Pareto sense, considering the points $(\| H\|,r(H))$.
 From the results in Table \ref{table3}, we observe:
\begin{itemize}
    \item[$\bullet$] On average, the number of solutions of each rank is much smaller for Penalized 1-norm. Both Cutting-plane and Cutting-plane Random generates the same or nearly that same
    number of solutions of each rank, which are on average, much greater than the numbers  for the other two methods.
    \item[$\bullet$] Nuclear-norm only obtains  solutions of minimum 1-norm for three ranks, and  for only one instance for each rank.
    \item[$\bullet$] Cutting-plane Random obtains the solution of minimum 1-norm in the vast majority of the cases, showing much superior results. The second best method is Cutting-plane.
    \item[$\bullet$] The results for the 0-norm are very similar to the results for the 1-norm, when analyzing Nuclear-norm and Penalized 1-norm, but when observing the the cutting-plane methods, we can see that the results are more balanced when comparing the 0-norms. 
\end{itemize}

		\begin{table}[!ht]
\centering{
	\begin{tabular}{c|cccc|cccc} %{c|rrrr|rrrr}
		\hline
		&\multicolumn{4}{c|}{$\|H\|_1$}&\multicolumn{4}{c}{$\|H\|_0$}\\
		Inst.&Nuclear-norm&Pen.1-norm&CP& CP.Rand&Nuclear-norm&Pen.1-norm&CP& CP.Rand\\
		\hline
		1  & 0 & 0 & 8 & 16 & 0 & 1 & 14 & 10 \\
2  & 0 & 0 & 5 & 15 & 0 & 1 & 7  & 11 \\
3  & 0 & 0 & 3 & 21 & 0 & 1 & 6  & 14 \\
4  & 0 & 4 & 2 & 17 & 0 & 2 & 5  & 14 \\
5  & 0 & 0 & 1 & 16 & 0 & 1 & 9  & 9  \\
6  & 0 & 1 & 4 & 12 & 0 & 1 & 6  & 9  \\
7  & 0 & 0 & 2 & 19 & 0 & 1 & 4  & 14 \\
8  & 1 & 0 & 9 & 7  & 1 & 1 & 5  & 8  \\
9  & 0 & 1 & 1 & 17 & 0 & 1 & 5  & 14 \\
10 & 1 & 1 & 2 & 16 & 0 & 1 & 8  & 10 \\
11 & 0 & 2 & 3 & 15 & 0 & 1 & 6  & 15 \\
12 & 0 & 0 & 1 & 20 & 0 & 1 & 11 & 15 \\
13 & 0 & 3 & 2 & 18 & 0 & 1 & 8  & 9  \\
14 & 0 & 0 & 1 & 18 & 0 & 1 & 6  & 14 \\
15 & 0 & 3 & 1 & 10 & 0 & 2 & 4  & 7  \\
16 & 0 & 0 & 6 & 16 & 0 & 1 & 7  & 17 \\
17 & 1 & 0 & 7 & 17 & 0 & 1 & 14 & 10 \\
18 & 0 & 0 & 3 & 21 & 0 & 1 & 8  & 15 \\
19 & 0 & 1 & 1 & 22 & 0 & 0 & 8  & 13 \\
20 & 0 & 0 & 3 & 19 & 0 & 1 & 4  & 18\\
		\hline
		\end{tabular}
			\caption{Number of solutions nondominated by other methods \label{table3a}}
		}
		\end{table}

Table \ref{table3a} presents results of the same experiments (as presented in Table \ref{table3}),
but now the rows indicate ``instance'' rather than ``rank'', and we tabulate the 
number of nondominated solutions generated by instance, in the Pareto sense, considering the points $(\| H\|,r(H))$.
We can clearly see that Cutting-plane Random is the best, with Cutting-plane also quite good.
Penalized 1-norm generates a few points that are nondominated, and Nuclear norm rarely does.

\begin{figure}[!ht]
	\centering
	\includegraphics[width=0.85\textwidth]{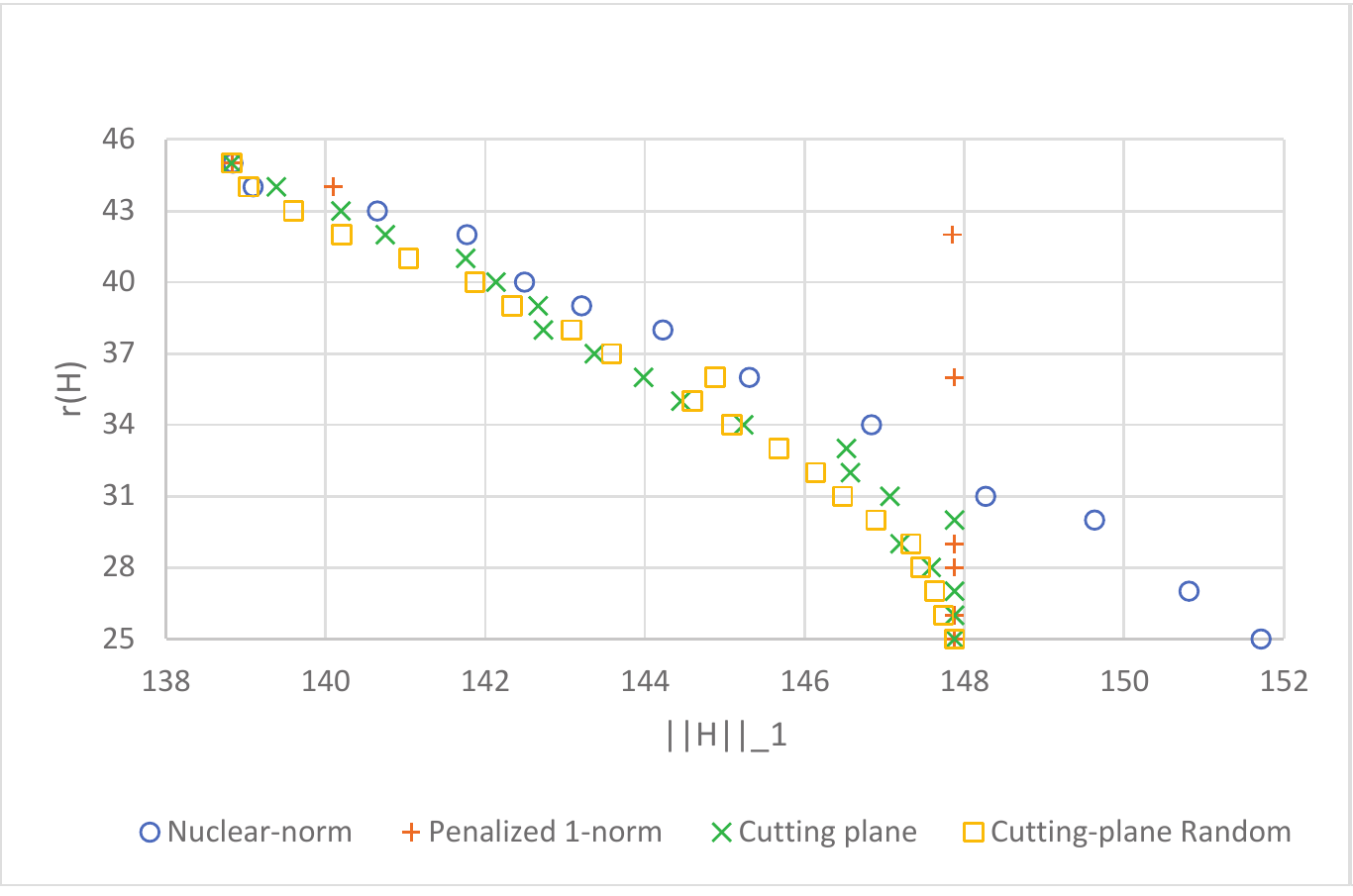}
	
	\medskip
	
	\includegraphics[width=0.85\textwidth]{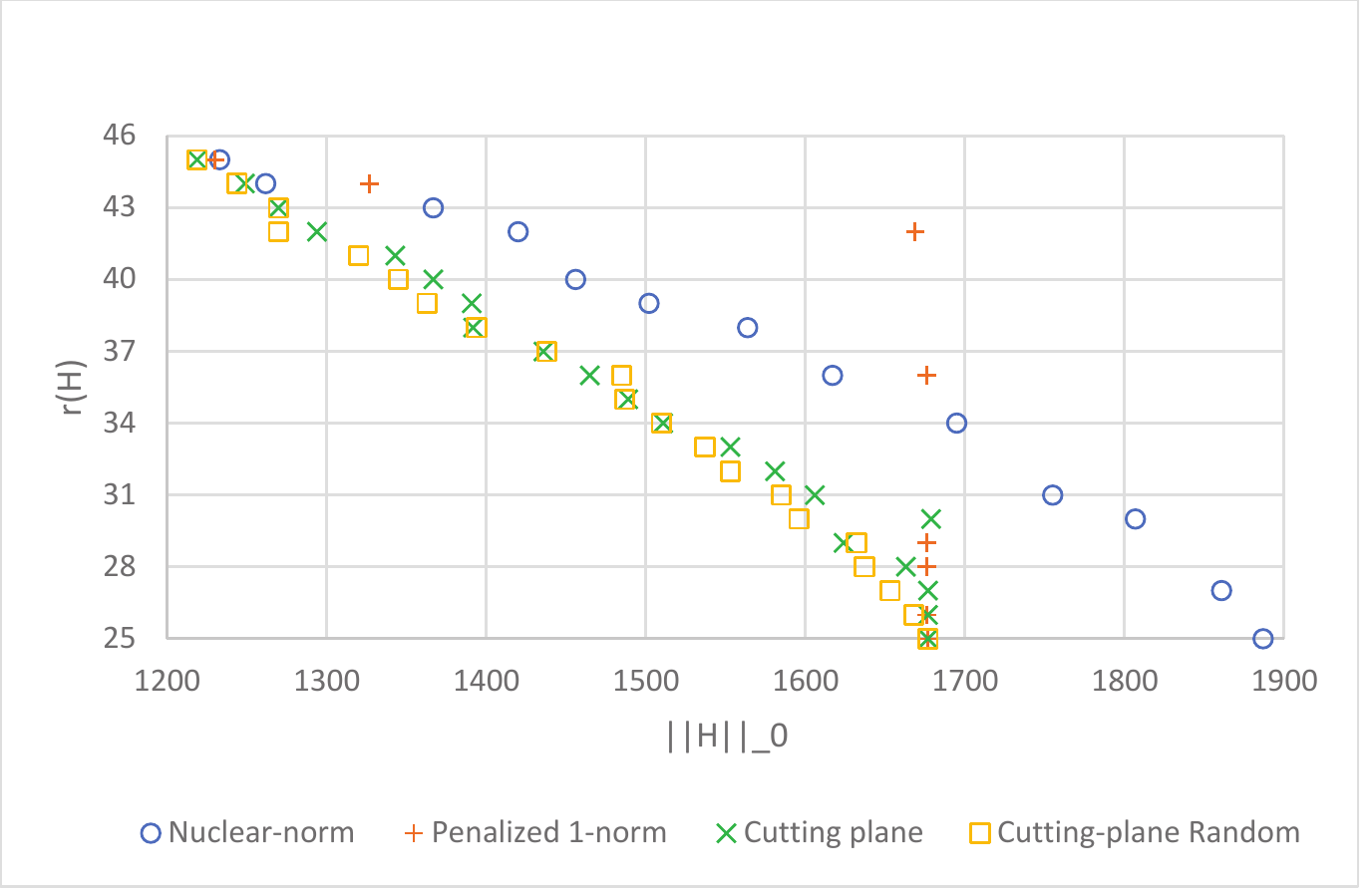}
	\caption{Pareto-curve approximations for an instance: four algorithms}\label{figure3}
\end{figure}

 The plots in Figure \ref{figure3} show the typical iterates generated by the four methods on one of our twenty random instances.
 These plots expand on the experiments summarized in ``Row 16'' of Table \ref{table3a}. 
We can see that Cutting-plane Random gives the best approximation of the Pareto curves, trading off 1/0-norms against rank.
Nuclear-norm performs quite poorly at generating $H$ with low 1/0-norms  for most ranks;
 it is particularly bad at the lower ranks. 
Even though Penalized 1-norm is overall quite poor, it can generate $H$ with decent 1/0-norms at very low and very high ranks.
Overall the Cutting-plane methods are the clear winners, with Cutting-plane Random offering some improvement over Cutting-plane.

\FloatBarrier

\section{Conclusion}\label{sec:conclusions}

% Comentar sobre a variacao de densidades e ranks. Vimos que mesmo com rank baixo r = 0.1m, se comporta de um jeito parecido r = 0.5m.
% Sobre as densidades: o rank apresenta mais variacoes (aumenta e diminui) com instancias esparsas. Alem disso, vi que em instancias mt esparsas tambem apresentam um aumento na esparsidade com o aumento na norma. O que nao e bom pra gente.
% \newline

There is a trade-off between 1-norm and sparsity against low-rank in the computation of ah-symmetric generalized inverses. These matrices can be applied in the solution of least-square problems, where  low 1-norm, sparsity, and low-rank are all desirable features in their numerical applications. In this work, we propose algorithmic approaches to numerically analyze this trade-off, using the fact that a generalized inverse has minimum rank if  and only if it satisfies the reflexive property. The algorithms start with a minimum 1-norm ah-symmetric generalized inverse and iteratively impose the reflexive property, 
gradually increasing  1-norm  and  decreasing  rank, until the minimum rank is obtained. The intermediate solutions from these algorithms can represent the best trade-off between 1-norm and sparsity against rank in numerical applications. 
Among the different strategies proposed to gradually impose the reflexive property, the best trade-off between norms and rank is obtained by a cutting-plane method that solves linear-programming problems at each iteration, 
applying a linear re-formulation of the reflexive property, obtained when combining it with the other properties of 
ah-symmetric generalized inverses.       
% We note that the numeric results used to construct the plots are shown in the spreadsheets  in our shared dropbox folder `IC Gabriel - P2 Relaxation/Spreadsheets'. 

% Six  groups of  plots are presented for the algorithms:
% \begin{enumerate}
% 		\item Cutting plane (most violated),
% 		\item Cutting plane (first violated),
% 		\item Augmented Lagrangian method,
% 		\item Penalized Frobenius,
% 	\item Penalized 1-norm.
% \end{enumerate}
% Each group shows average results for the 5 instances with the same configuration, as described above.
% These instances were not used in experiments with the algorithms:
% \begin{enumerate}
% 	\setcounter{enumi}{5}
% 	\item Lagrangian method,
% 	\item Penalized rank.
% \end{enumerate}

 \FloatBarrier

\bibliographystyle{plain+eid}


%%%%%%\bibliographystyle{spmpsci}      % mathematics and physical sciences

\bibliography{ginv}

%\newpage
%\newpage
%\newpage
%
%\newpage
%\bibliographystyle{amsplain}

%\bibliography{ginv}

\end{document}